\documentclass[fleqn,a4paper]{article}

\begin{document}
\raggedright
\sffamily

\title{On the repeated inversion of a covariance matrix}
\author{M.\ de Jong\\
  {\small NWO-I, Nikhef, PO Box 41882, Amsterdam, 1098 DB Netherlands} \\
  {\small Leiden University, Leiden Institute of Physics, PO Box 9504, Leiden, 2300 RA Netherlands}
}
\maketitle

\begin{abstract}
  In many cases, the values of some model parameters are determined 
  by maximising the likelihood of a set of data points given the parameter values. 
  The presence of outliers in the data and correlations between data points
  complicate this procedure.
  An efficient procedure for the elimination of outliers is presented
  which takes the correlations between data points into account. 
\end{abstract}

\section{Introduction}

In general, the values of some model parameters are optimally determined 
by maximising the likelihood of a set of data points given the parameter values
(see e.g.\ \cite{ref:statistical-methods}).
In this, a data point that has to low a probability to match the model can be considered an outlier.
The presence of such outliers in the data can readily be accommodated in
the probability density function of the data points.
It is, however, not straight forward to also take the correlations between the data points into account. 
If the underlying probability density functions are normal distributions, 
both the uncorrelated and correlated uncertainties of the data points 
can be incorporated in a single matrix.
The values of the model parameters can then be determined by minimising the $\chi^2$:

\begin{eqnarray}
  \chi^2 & = & \epsilon^{T} \; V^{-1} \; \epsilon
  \label{eq:chi2}
\end{eqnarray}

where
$\epsilon$ and $\epsilon^T$ are the $(N \times 1)$ and $(1 \times N)$ vectors
containing the distances between the model and the data points and
$V$ is a $(N \times N)$ matrix.
The elements of $V$ are set as follows.

\begin{eqnarray}
  V_{ii} & = & (\sigma_i)^2                                \\
    \label{eq:matrix-on}
  V_{ij} & = & \sum_{k}
  \frac{\partial \epsilon_i}{\partial u_k}
  \frac{\partial \epsilon_j}{\partial u_k}(\delta u_k)^2
  \label{eq:matrix-off}
\end{eqnarray}

where
$\sigma_i$ refers to the uncorrelated uncertainty of data point $i$ and 
$u_k$ to some correlation parameter which itself has an an uncertainty $\delta u_k$
The terms in the summation of equation \ref{eq:matrix-off} are commonly referred
to as the covariances of the data points.
The matrix $V$ is therefore often called the covariance matrix.
By construction, the matrix $V$ is symmetric and can thus be inverted using an LDU decomposition.
The computation of $V^{-1}$ then requires $\mathcal{O}(N^3)$ operations. \\~\\

The presence of outliers cannot easily be incorporated in the covariance matrix.
So, it is desirable to remove the outliers and repeat the fit.
For this, a criterion is required to identify an outlier.
A common criterion is based on the value of the so-called standard deviation, $D$:

\begin{eqnarray}
  D_k  & \equiv &  \frac{|\epsilon_k|}{\sigma_k}
  \label{eq:outlier-1}
\end{eqnarray}

A typical maximal allowed value of $D$ is $D_{\mathrm{max}} \; = \; 3$
which --in the absence of correlations-- corresponds to a probability to keep a good data point of about $0.997$.
The removal of outliers could simply proceed by
removing the data point with the largest $D$ and repeating the fit 
until there are no more data points with $D > D_{\mathrm{max}}$. \\~\\

In case the correlations between the data points are strong,
this procedure may no longer be adequate.
In this scenario, the absolute values of some off-diagonal elements of $V$
are comparable to (or even larger than) the values of the diagonal elements.
The standard deviation is then no longer a good criterion because 
the distances and covariances of the other data points should also be taken into account.
A brute force procedure to identify outliers is to
\emph{1)} remove data point $k$,
\emph{2)} determine $V'$,
\emph{3)} invert $V'$ and
\emph{4)} minimise $\chi^2$, 
for each data point $k$.
In this, $V'$ corresponds to the $(N-1)\times(N-1)$ covariance matrix.
The change in $\chi^2$ is then a good criterion to identify an outlier, i.e: 

\begin{eqnarray}
    D_k  & \equiv &  \sqrt{\chi^2 - \chi^2_k}
  \label{eq:outlier-2}
\end{eqnarray}

where $\chi^2_k$ corresponds to the $\chi^2$ of the fit after removal of data point $k$.
As each inversion of $V'$ requires $\mathcal{O}((N-1)^3)$ operations, 
this way of eliminating outliers requires $\mathcal{O}(N^4)$ operations.
For a large number of data points, the number of operations
needed to eliminate outliers may become too excessive.\\~\\

An alternative exist, based on the fact that the removal of data point $k$
is equivalent to setting the corresponding uncertainty $\sigma_k$ to infinity.
By doing so, the $(N \times N)$ covariance matrix $V'$ can be decomposed as follows:

\begin{eqnarray}
  V' & = &  V  ~+~  g \; \delta_{k,k}                                         
\end{eqnarray}

In this, 
$g$ is some arbitrary large value; in any case much larger than $\sigma_k$.
The matrix $\delta_{i,j}$ has $1$ at row $i$ and column $j$ and $0$ everywhere else.
To repeat the fit without data point $k$, it is required to invert matrix $V'$.
As a first step, the known inverse of the original matrix $V$ is considered.
This is possible because $V^{-1}$ and $V'$ have the same dimensions. 

\begin{eqnarray}
   V' \times V^{-1} &   =    &  (V + g \delta_{k,k}) \times V^{-1}                \\
   \mbox{}                                                                       \nonumber\\
   &   =    &  VV^{-1} + g \; \delta_{k,k} V^{-1}                                  \\
   \mbox{}                                                                       \nonumber\\
   &   =    &  I ~+~  g 
   \left( \begin{array}{ccccc}
     0           & \multicolumn{3}{c}{\cdots}   &       0      \\
     \vdots      &        &             &       &    \vdots    \\
     0           & \multicolumn{3}{c}{\cdots}   &       0      \\
     V^{-1}_{k,1} & \cdots & V^{-1}_{k,k} & \cdots  & V^{-1}_{k,N} \makebox[0cm][l]{\hspace*{4\arraycolsep}$\leftarrow \textrm{row}~k$}\\
     0           & \multicolumn{3}{c}{\cdots}   &       0      \\
     \vdots      &        &             &       &    \vdots    \\
     0           & \multicolumn{3}{c}{\cdots}   &       0      \\
   \end{array} \right)                                                            \\
   \mbox{}                                                                        \nonumber\\
   & \equiv &  A                                   
\end{eqnarray}

where $I$ is the identity matrix.
It is obvious that the inverse of matrix $V'$ is equal to the product $V^{-1} A^{-1}$.
Consequently, the problem of inverting $V'$ is reduced to inverting matrix $A$.
The inverse of matrix $A$ is trivial, namely:

\begin{eqnarray}
  A^{-1}        &   =    &  I ~-~  \frac{g}{1 + g V^{-1}_{k,k}}
  \left( \begin{array}{ccccc}
    0           &  \multicolumn{3}{c}{\cdots}   &       0      \\
    \vdots      &        &             &        &    \vdots    \\
    0           &  \multicolumn{3}{c}{\cdots}   &       0      \\
    V^{-1}_{k,1} & \cdots & V^{-1}_{k,k} & \cdots & V^{-1}_{k,N} \makebox[0cm][l]{\hspace*{4\arraycolsep}$\leftarrow \textrm{row}~k$}\\
    0           &  \multicolumn{3}{c}{\cdots}   &       0      \\
    \vdots      &        &             &        &    \vdots    \\
    0           &  \multicolumn{3}{c}{\cdots}   &       0      \\
  \end{array} \right)
  \label{eq:A-1}
\end{eqnarray}

Now, one can let go $g \rightarrow \infty$ and multiply $V^{-1}$ and $A^{-1}$
to obtain the inverse of matrix $V'$.
As can be seem from equation \ref{eq:A-1},
the product of $V^{-1}$ and $A^{-1}$ requires $\mathcal{O}(N^2)$ operations.
Hence, this way of eliminating outliers only requires $\mathcal{O}(N^3)$ operations
which actually is equal to the number of operations needed to invert the covariance matrix for the first fit.
It is interesting to note that no additional memory is required for the computation of $V'^{-1}$.\\~\\

A further reduction in the number of operations is possible when the outcome of the first fit is retained.
In that case, the distance $\epsilon_i$ between the model and data point $i$ stays the same.
The $\chi^2$ without data point $k$ can then be expressed as:

\begin{eqnarray}
 \chi^2_{k}  & = &  \sum_{i=1}^{N} \sum_{j=1}^{N} ~ \epsilon_i ~ V'^{-1}_{i,j} ~ \epsilon_j                   \\
 \mbox{}                                                                                                   \nonumber\\
 \mbox{}     & = &  \sum_{i=1}^{N} \sum_{j=1}^{N} ~ \epsilon_i ~ (V^{-1}A^{-1})_{i,j} ~ \epsilon_j            \\
 \mbox{}                                                                                                   \nonumber\\
 \mbox{}     & = &  \chi^2 - \frac{1}{V^{-1}_{k,k}} \left(\sum_{j=1}^{n} V^{-1}_{k,j} ~ \epsilon_j\right)^2
 \label{eq:chi2-k}
\end{eqnarray}

As can be seen from equation \ref{eq:chi2-k},
the distances and covariances of the other data points are taken into account
but the inverse of $V'$ is no longer required.
As a result, the summation takes only $\mathcal{O}(N)$ operations and
the elimination of outliers $\mathcal{O}(N^2)$ operations.

\section{Conclusions}

An efficient method is presented to eliminate outliers from a set of data points 
which takes the correlations between the data points into account.
The number of operations needed for this procedure is 
the same as that needed for the inversion of the covariance matrix.

\bibliographystyle{plain}
\bibliography{Matrix}

\begin{thebibliography}{1}

\bibitem{ref:statistical-methods}
{W.J.~Metzger}.
\newblock {\em \mbox{STATISTICAL METHODS IN DATA ANALYSIS}}.
\newblock {\\*[0cm]Katholieke Universiteit Nijmegen, the Netherlands}, 1996.
\newblock {HEN-343}.

\end{thebibliography}

\end{document}